# Intelligent Home Energy Management System for Distributed Renewable Generators, Dispatchable Residential Loads and Distributed Energy Storage Devices

Adetokunbo Ajao, Jingwei Luo, Zheming Liang
Electrical and Computer Engineering
University of Michigan-Dearborn, USA

Qais H. Alsafasfeh
Electrical Engineering
Islamic Educational, Scientific
And Cultural Organization (ISESCO)

Wencong Su
Electrical and Computer Engineering
University of Michigan-Dearborn, USA
Email: wencong@umich.edu

*Abstract*—This paper presents an intelligent home energy management system integrated with dispatchable loads (e.g., clothes washers and dryers), distributed renewable generators (e.g., roof-top solar panels), and distributed energy storage devices (e.g., plug-in electric vehicles). The overall goal is to reduce the total operating costs and the carbon emissions for a future residential house, while satisfying the end-users' comfort levels. This paper models a wide variety of home appliances and formulates the economic operation problem using mixed integer linear programming. Case studies are performed to validate and demonstrate the effectiveness of the proposed solution algorithm. Simulation results also show the positive impact of dispatchable loads, distributed renewable generators, and distributed energy storage devices on a future residential house.

I. INTRODUCTION

Commercial and residential buildings consumed almost 40% of the primary energy and approximately 70% of the electricity in the United States in 2012, and the trend continues to escalate. Intelligent home energy management is a viable solution to reduce energy costs, maintain customer comfort levels, accommodate the integration of distributed renewable energy resources, and facilitate demand-side management (DSM) and demand response (DR) programs [1]-[3].

The emerging technology of distributed energy storage devices (DESD), such as plug-in electric vehicles (PHEVs) [4], opens new possibilities for bi-directional power flows. For example, vehicle-to-grid (V2G) and vehicle-to-home (V2H) allow a PHEV battery to feed energy directly back into the power grids and/or power other local loads. In addition, the introduction of dynamic electricity pricing further helps customers to reduce the cost of the electricity they use [5]. Recently, common dynamic pricing schemes include time-of-use pricing (TOU), critical peak pricing, and real-time pricing (RTP).

In [6], the authors proposed a scheduling method for household appliances based on dynamic price signals. In [7] the authors proposed different types of household appliance models, but did not explore the bi-directional power flow in which the consumer can sell electricity back to the power grids. The authors in [8] proposed a bi-directional PHEV charging/discharging model, but did not analyze the effect of household appliances on the consumption of electricity.

The above-mentioned literature survey claims an urgent need for an intelligent home energy management system integrated with dispatchable loads (e.g., clothes washers and dryers), distributed renewable generators (e.g., roof-top solar panels), and DESDs (e.g., PHEVs). Fig. 1 shows the envisioned architecture and major components of a residential house.

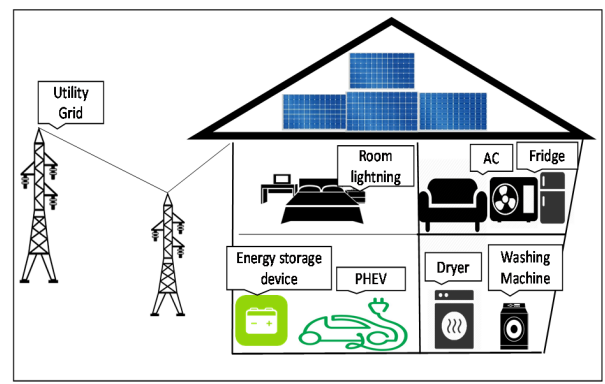

Fig. 1. Envisioned architecture of a future residential house.

The major contributions of this paper can be summarized as follows:

- Modeling of a wide variety of household appliances and their operating constraints;
- Formulation of an objective function to minimize the operation cost considering bi-directional power flow and customer preferences;
- Analysis of the impact of dispatchable loads, distributed renewable generators, and distributed energy storage devices on a future residential house.

The rest of the paper is organized as follows: Section II discusses the mathematical modeling and the objective function formulation. Section III describes the system configuration and the input data. Section IV presents case studies and investigates the impact of dispatchable loads, distributed renewable generators, and DESDs on a future residential house. Section V highlights the conclusion and future work.

## II. PROBLEM FORMULATION

### A. Objective Function

The primary goal of this research is to minimize the total operating cost of a future residential house. The objective function is as follows:

$$Cost_{house} = \sum_{t=1}^{24}[(P_t^b Pr_t^b - P_t^s Pr_t^s)] + \sum_{t=1}^{24}[(P_t^{CH,DESD}\eta_1 + P_t^{DIS,DESD}\eta_1) + (P_t^{CH,PHEV}\eta_2 + P_t^{DIS,PHEV}\eta_2)] \quad (1a)$$

$$Cost_{DESD\&PHEV} = \sum_{t=1}^{24}[((P_t^{CH,DESD}\eta_1 + P_t^{DIS,DESD}\eta_1) + (P_t^{CH,PHEV}\eta_2 + P_t^{DIS,PHEV}\eta_2))] \quad (1b)$$

Equation (1a) shows the total operating cost of the residential green house, $P_t^b$ (kW) denotes the power bought from the power grid at time interval t, while $Pr_t^b$ denotes the dynamic electricity price ($/kWh) offered by the utility grid, $P_t^s$ (kW) is the amount of power sold from the PV, DESD, and the PHEV at time interval t, and $Pr_t^s$ (kWh) is a flat price agreed to by the consumer at which they will sell their power back to the grid. The bi-directional flow of PHEVs and DESDs makes it possible for a consumer to meet their household demand and still have excess energy at a time interval to sell back to the grid. Equation (1b) shows the degradation cost of the DESD and the PHEV ($Cost_{DESD\&PHEV}$), which is part of the total cost of operating the house, $P_t^{CH,DESD}$ and $P_t^{DIS,DESD}$ (kW) indicate the charging power and discharging power of the DESD respectively, while $P_t^{CH,PHEV}$ and $P_t^{DIS,PHEV}$ (kW) indicate the charging and the discharging power of the PHEV respectively, $\eta_1$ ($/kwh) is the degradation cost associated with the charging and discharging process of the DESD, while $\eta_2$ ($/kwh) is the degradation cost associated with the charging and discharging process of the PHEV.

### B. System Constraints

Load in our model can be divided into controllable and uncontrollable load. Controllable load is the load that the operation schedule can shift or manipulate without affecting the operation itself, which has been discussed elaborately in [9], while uncontrollable load is the load that cannot be shifted because it will affect the mode of operation. In this paper the house contains an AC (air conditioner), fridge, washer and dryer as the basic load, in which the washer and dryer are categorized as controllable load and the fridge and AC are uncontrollable loads.

*1) Air conditioner and heater operating constraints*

The AC and heater work on the principle in which when the inside temperature is higher than a specific temperature, the AC tends to switch ON and if it is lower than a specified temperature the heater switches ON. If the inside temperature of the house, $\theta_{in}$, is at an optimal level based on the customer preference, neither of the devices will run. $\theta_{in}$ has an effect on the power consumed by the AC and heater, which can be determined by:

$$\theta_{in} = \theta_{in}(t-1) + [v_{ac} A(t) - u_{ac} s_{ac(t)} + I_{ac}(\theta_{out}(t) - \theta_{in}(t-1))], \forall t \in T \quad (2)$$

$$\theta_{in} = \theta_{in}(t-1) + [v_{ht} A(t) - u_{ht} s_{ht}(t) + I_{ht}(\theta_{out}(t) - \theta_{in}(t-1))], \forall t \in T \quad (3)$$

$$S_{ac}(t) + S_{ht}(t) \leq 1, \forall t \in T \quad (4)$$

In (2) and (3) $\theta_{in}$ at any time interval depends on the previous temperature and the cooling and warming effect of the OFF state of the AC $v_{ac}$ and the heater $v_{ht}$. A(t) is the activity level of the house (%) for both equations which is essentially the energy consumption of a typical household. $u_{ac}$ and $u_{ht}$ is the cooling and warming effect of the ON state of the AC and heater respectively, $I_{ac}$, and $I_{ht}$, is the effect of outside and inside temperature difference on the AC and heater respectively which has a huge impact on the inside temperature per time interval. (4) indicates that the AC and the heater are not to operate at the same time, where $S_{ac}$ and $S_{ht}$ represent the ON/OFF status of the AC and the heater respectively, which can only take binary values.

*2) Fridge operation constraints*

$$\theta_{fr} = \theta_{fr}(t-1) + [u_{fr} A_{fr}(t) - v_{fr} s_{fr}(t) + \alpha_{fr}], \forall t \in T \quad (5)$$

In equation (5), $\theta_{fr}$ indicates the inside temperature of the fridge at a time interval which is dependent on the previous temperature of the fridge $\theta_{fr}(t-1)$, the effect of the activity level on the fridge temperature $u_{fr}$, in which as the activity of the house increases there will be more demand on cooling of the fridge, $A_{fr}(t)$ indicates the activity level of the fridge, while $v_{fr}$, indicates the cooling and warming effect of the ON state of the fridge and $s_{fr(t)}$ is the binary variable that controls the switching of the fridge, moreover $\alpha_{fr}$, indicates the warming effect of the OFF state of the fridge. $\theta_{fr}$ denotes the effect of the fridge on power consumption, in which the higher the inside temperature of the fridge the higher the power consumption.

*3) Washer and dryer operation constraints*

As discussed earlier, the washer and dryer are categorized as a controllable load, where i represent either the washer or the dryer. The operating constraints of the washer and dryer are shown below:

$$u_i(t) - v_i(t) = s_i(t) - s_i(t-1), \forall t \in T, \forall i \in I \quad (6)$$

$$u_i(t) + v_i(t) \leq 1, \forall t \in T, \forall i \in I \quad (7)$$

$$\sum_{t \in T_i} s_i(k) = O_i^{rt}, \forall t \in T \quad (8)$$

$$\sum_{k=t}^{t+O_i^{mst}} s_i(k) \leq O_i^{mst} + M(1-u_i(t)), \forall t \in T \quad (9)$$

$$\sum_{k=t-U_i+1}^{t} u_i(t) \leq s_i(t), \forall t \in T \quad (10)$$

$$\sum_{k=t-D_i+1}^{t} u_i(t) \leq 1 - s_i(t), \forall t \in T \quad (11)$$

$$s_{dryer}(t) \leq \sum_{k=1}^{\Omega} s_{washer}(t-k), \forall t \in T \quad (12)$$

$$s_{dryer}(t) + s_{washer}(t) \leq 1, \forall t \in T \quad (13)$$

Equations (6) and (7) model the shut down and start up constraints for the washer and dryer in order not to damage the device, where $u_i(t)$ and $v_i(t)$ represent the binary variable denoting the startup of the device i at time t and the shutdown of the device i at time t respectively. Equation (8) constrains the device to operate at a particular operation time, $O_i^{rt}$. Equation (9) constrains the device to operate at a maximum successive time, $O_i^{mst}$, where M represents a large positive number. Equations (10) and (11) constrain the device to operate at a minimum up time, $U_i$, and minimum down time, $D_i$, respectively. Moreover, (12) sets the dryer to operate after the washer has finished its task, and $\Omega$ represents the time gap as set by the customer preference. Equation (13) constrains the dryer and washer not to operate at the same time.

*4) Power balance equations*

$$P_t^b + P_t^{used,PV} + P_t^{used,PHEV} + P_t^{used,DESD} = P_t^{RGH} + P_t^{CH,PHEV} + P_t^{CH,DESD} + P_t^s, \forall t \in T \quad (14)$$

Equation (14) shows the power balance equation, where $P_t^b$ and $P_t^s$ represent the total power bought and sold from or to the grid at a time interval. $P_t^{used,Pv}, P_t^{used,PHEV}$ and $P_t^{used,DESD}$ represent the power supplied to meet the residential load demand at a time interval by the PV, PHEV and DESD respectively. $P_t^{RGH}$ represents the total power consumed by the fridge, AC, washer and dryer at a time interval. It should be noted that $P_t^{CH,PHEV}$ and $P_t^{CH,DESD}$ are the charging power demands of the PHEV and DESD, respectively, in which these devices operate as an electrical generator when discharging and energy consuming load when charging.

*5) DESD operation constraints*

The bi-directional flow constraints of the DESD are shown below. The DESD charges when the dynamic electricity price is low and discharges to meet the household load when the electricity price is low.

$$P_t^{used,DESD} + P_t^{s,DESD} = P_t^{DIS,DESD} \eta_{DD}, \forall t \quad (15a)$$

$$P_t^{CH,DESD} \leq CR_t^{DESD} s_t^{DESD}, \forall t \quad (15b)$$

$$P_t^{DIS,DESD} \leq DR_t^{DESD}(1 - s_t^{DESD}), \forall t \quad (15c)$$

$$SOE_t^{DESD} = SOE_{t-1}^{DESD} + P_t^{CH,DESD} \eta_{CD} - P_t^{DIS,DESD}, \forall t \quad (15d)$$

$$SOE_1^{DESD} = SOE^{DESD,INI}, \forall t \quad (15e)$$

$$SOE_t^{DESD} \leq SOE^{DESD,MAX}, \forall t \quad (15f)$$

$$SOE_t^{DESD} \geq SOE_t^{DESD,MIN}, \forall t \quad (15g)$$

Equation (15a) indicates the DESD can be used to meet the household load demand where $P_t^{s,DESD}$ represents the power sold by the DESD at a time interval t, $P_t^{DIS,DESD}$ represents the discharging power, and $\eta_{DD}$ represents the discharging efficiency of the DESD. (15b) helps to limit the charging power of the DESD, where $P_t^{CH,DESD}$ is the charging power, $CR_t^{desd}$ is the charging rate, and $s_t^{DESD}$ is the ON/OFF status of the DESD. This constraint protects the life cycle of the DESD. (15c) controls the discharging rate of the DESD, where $P_t^{DIS,DESD}$ represents the discharging power of the DESD and $DR_t^{DESD}$ represents the discharging rate of the DESD at time interval t. (15d) represents the state of energy of the DESD, which is dependent on the previous SOE, the charging efficiency, $\eta_{CD}$, and the charging power deducted from the discharged power at a time interval t. (15e) makes the SOE at the first interval equal to the initial state of energy, $SOE^{DESD,INI}$. Equations (15f) and (15g) constrain the SOE of the DESD to not exceed a maximum $SOE^{DESD,MAX}$ and not go below a minimum $SOE^{DESD,MIN}$ respectively.

*6) PHEV operation constraints*

$$P_t^{used,PHEV} + P_t^{s,PHEV} = P_t^{DIS,PHEV} \eta_{DP}, \forall t \in [T^a, T^d] \quad (16a)$$

$$P_t^{CH,PHEV} \leq CR_t^{PHEV} s_t^{PHEV} \forall t \in [T^a, T^d] \quad (16b)$$

$$P_t^{DIS,PHEV} \leq DR_t^{PHEV} (1 - s_t^{PHEV}), \forall t \in [T^a, T^d] \quad (16c)$$

$$SOE_t^{PHEV} = SOE_{t-1}^{PHEV} + P_t^{CH,PHEV} \eta_{CP} + P_t^{DIS,PHEV}, \forall t \in [T^a, T^d] \quad (16d)$$

$$SOE_{T^a}^{PHEV} = SOE^{PHEV,INI} \quad (16e)$$

$$SOE_t^{PHEV} \leq SOE_t^{PHEV,MAX}, \forall t \in [T^a, T^d] \quad (16f)$$

$$SOE_t^{PHEV} \geq SOE_t^{PHEV,MIN}, \forall t \in [T^a, T^d] \quad (16g)$$

$$SOE_{T^{f,c}}^{PHEV} = SOE_t^{PHEV,MAX} \quad (16h)$$

The mode of operation of the PHEV is constrained to operate at particular $T^a$ and $T^d$, which correspond to the arrival and departure time of the PHEV, respectively. Equation (16a) indicates the PHEV can be used to meet the residential green house (RGH) load demand, where $P_t^{s,PHEV}$ represents the power sold by the PHEV at a time interval t, $P_t^{DIS,PHEV}$ represents the discharging power, and $\eta_{DP}$ represents the discharging efficiency of the PHEV. (16b) helps to limit the charging power of the PHEV, where $P_t^{CH,PHEV}$ is the charging power, $CR_t^{PHEV}$ is the charging rate, and $s_t^{PHEV}$ is the ON/OFF status of the PHEV. This constraint protects the life cycle of the PHEV. (16c) controls the discharging rate of the PHEV, where $P_t^{DIS,PHEV}$ represents the discharging power of the PHEV and $DR_t^{PHEV}$ represents the discharging rate of the PHEV at time interval t. (16d) represents the state of energy of the PHEV, which is dependent on the previous SOE, the charging efficiency, $\eta_{CP}$, and the charging power deducted from the discharged power at a time interval t. (16e) makes the SOE at arrival time $T^a$ of the driver equal to the initial $SOE^{PHEV,INI}$. (16f) and (16g) constrain the SOE of the EV not to exceed a maximum $SOE^{PHEV,MAX}$ and not fall below a minimum $SOE^{PHEV,MIN}$, respectively, within $T^a$ and $T^d$ of the driver while connected to the bi-directional smart meter. (16h) indicates the period, $T^{f,c}$, in which PHEV must be fully charged. It should also be noted that the PHEV constraints won't be initiated until the arrival time and departure time of the driver.

*7) PV operation constraints*

$$P_t^{used,PV} + P_t^{s,PV} + P_t^{g,PV}, \forall t \quad (17)$$

Equation (17) shows that the PV can be used in meeting the RGH load demand, where $P_t^{g,PV}$ represents the total generation power of the solar panel

*8) Total power injected into the grid*

$$P_t^s = P_t^{s,PV} + P_t^{s,DESD} + P_t^{s,PHEV}, \forall t \quad (18)$$

Equation (18) shows the total power sold, $P_t^s$, by the PV, DESD and PHEV respectively.

*9) Power transaction regulation*

In order to control the amount of power that can be sold and received from the grid, a set of constraints was modeled as shown below:

$$P_t^b \leq m1 s_t^{grid}, \forall t \quad (19)$$

$$P_t^s \leq m2(1 - s_t^{grid}), \forall t \quad (20)$$

Where $m1$ and $m2$ show the maximum amount of power that can be received from the grid and the maximum amount of power that can be sold to the grid, respectively. These constraints help to keep the bi-directional flow of power within the customer preference and at an agreed rate between the grid and the consumer. $s_t^{grid}$ denotes the ON/OFF status of the grid, which will be 1 if the RGH is taking power from the grid, or 0 if the RGH is not taking power from the grid or is in an idle state.

III. NUMERICAL SETTINGS

*A. Temperature*

A particular set of summer outside temperature data was taken from Atlanta, Georgia in the United States [10]. Fig. 2 shows the distribution of the outside temperature, which ranges from 23-29 degrees Celsius.

*B. Activity Level*

The activity level of the RGH is determined by calculating the hourly energy consumption of a typical resident in Atlanta, Georgia in terms of a percentage. Fig. 3 shows the activity level of a residential house in Georgia, in which there is an increase of activity level at 2am and 11am, and a clear reduction during the night, as the occupant has less power demand [11].

*C. Dynamic Electricity Price*

In order to test the variability of the effect of PHEVs, DESDs and PVs on energy cost, a dynamic price is implemented. Fig. 4 shows the dynamic electric price [12].

*D. Solar Power Generation*

A 1-kW solar panel, as used in Atlanta during a typical summer, is used in this analysis. The highest output was at 1pm, at 0.82 kW. Fig. 5 shows the hourly solar power generation in Atlanta [13].

TABLE I. PARAMETER VALUES OF HOUSING APPLIANCES [14]

| Device | Parameter |
|---|---|
| AC | Rated power is 1.9kw |
| Fridge | Rated power is 0.42kw |
| Washer | Rated power is 0.5kw |
| Dryer | Rated power is 3.5kw |

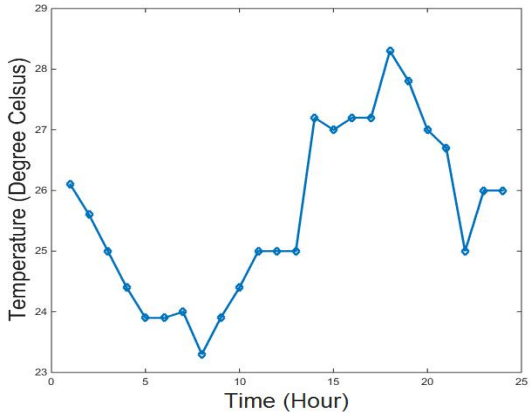

Fig. 2. Outside temperature in Atlanta, Georgia.

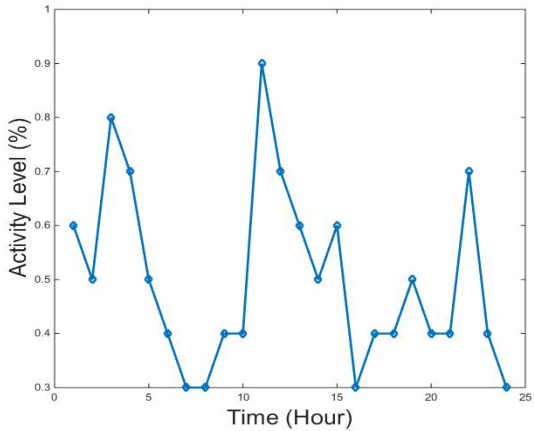

Fig. 3. Activity level of a typical residential house in Georgia.

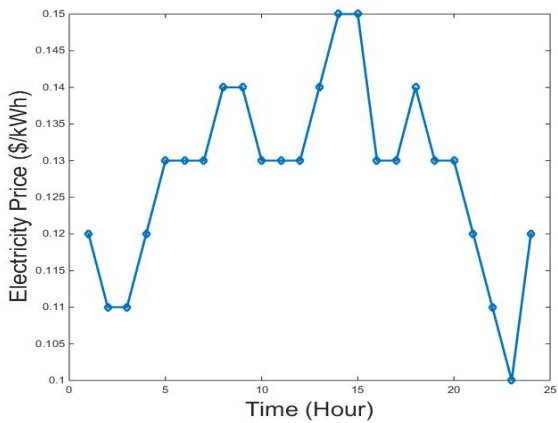

Fig. 4. Dynamic electricity price.

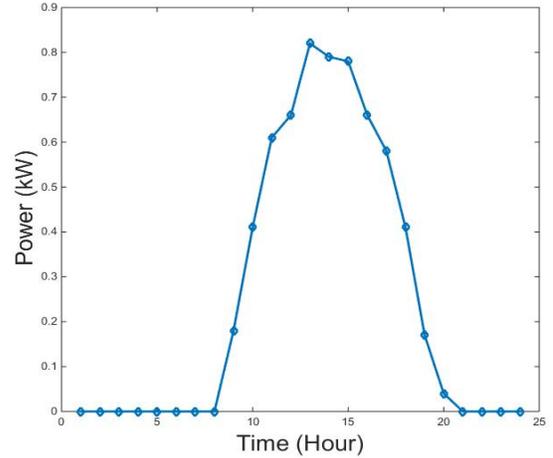

Fig. 5. Solar power generation in Atlanta (1-kW generator).

## IV. CASE STUDIES

The formulation of the objective function and the system constraints from (1-19) was implemented in Matlab and solved by the IBM ILOG CPLEX optimizer [15]. All the simulations were run on a computer with an Intel core i7-4790 CPU@ 3.60GHz with a 16.00 GB memory.

### A. Scenario I

In order to demonstrate the effect of PHEV and DESD bidirectional flow, which helps to reduce the total operating cost, a scenario in which no PV or DESD is initially used consisting only of a household basic load, which consists of the fridge, washer and dryer, heater and AC was considered. A 1kw solar panel was used and a 1kw energy storage device (DESD) was used as well. The driver is assumed to be back at 6pm and is willing to charge as soon as the dynamic electricity price is low. A Chevy Volt of 16kwh [12] is used, which has both a charging and discharging rate of 3.3kwh.

### B. Scenario II

In this scenario, PV, DESD and PHEV bi-directional flow are considered to show the reduction in total operating cost for the day. The PV, DESD and PHEV can sell power to the grid, and the DESD and PHEV can receive power from the grid to charge their batteries. The same household load was considered in the scenario as well. In this scenario, the driver of the car is assumed to arrive at 6pm and discharge first, charging later when the electricity price is low. The next section discusses the results of the simulation of the two scenarios.

## V. SIMULATION RESULTS AND DISCUSSION

After running the simulation, the total operating cost of the RGH for the day was $4.98 for scenario I and $ 4.09 for the scenario II, which shows a difference of $0.89. Fig. 5 shows the state of energy of the PHEV for the two scenarios. It can be seen in Fig. 6 that in scenario I, the PHEV kept charging despite the high dynamic electricity price, while in scenario II, the PHEV discharged some of its energy to meet the

household load before later charging when the dynamic price was low. This factor contributed to the total operating cost of Scenario I. Furthermore, Fig. 7 shows the total power bought by the RGH for the two scenarios. It can be seen from the figure that more power was bought from the grid in scenario I during the peak period than in Scenario II, in which the combination of the bidirectional flow of the DESD and PHEV with the PV helped to reduce the amount of power bought from the grid. Table II shows the total cost of operating the RGH for a 24-hr interval. It can be seen that it costs less to operate the RGH when taking advantage of the effect of bidirectional ESD and PHEV.

## VI. CONCLUSION

This paper proposes the use of PHEV and DESD bidirectional flow to help reduce the cost of operating the RGH, in which the PHEV, DESD and PV can sell energy to the grid, while only the PHEV and DESD can receive energy from the grid to charge their batteries. Consensus-based distributed control methods will be our future research work [17] because they have received more attention recently [18].

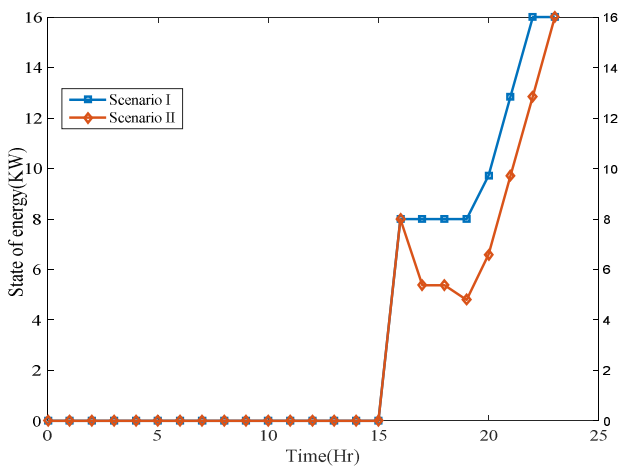

Fig. 6. State of energy of the PHEV in the two scenarios.

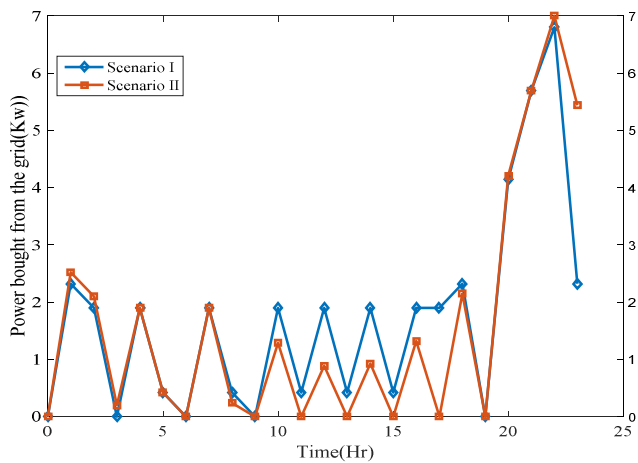

Fig. 7. Power bought from the grid in the two scenarios.

TABLE II. TOTAL COST OF DAILY OPERATION

| No PV, DESD | PV, DESD and PHEV bidirectional flow |
|---|---|
| $4.98 | $ 4.09 |